\documentclass{article}
\usepackage{amsmath,amssymb}
\newcommand{\Z}{\ensuremath{\mathbf Z}}
\newcommand{\N}{\ensuremath{ \mathbf N }}
 
\newtheorem{theorem}{Theorem}
\newcommand{\bt}{\begin{theorem}}
\newcommand{\et}{\end{theorem}}
\newtheorem{lemma}{Lemma}
\newcommand{\bl}{\begin{lemma}}
\newcommand{\el}{\end{lemma}}
\newcommand{\pf}{{\bf Proof}.\ }
\newcommand{\be}{\begin{eqnarray}}
\newcommand{\ee}{\end{eqnarray}}
\newcommand{\beq}{\begin{equation}}
\newcommand{\eeq}{\end{equation}}
\newcommand{\benum}{\begin{enumerate}}
\newcommand{\eenum}{\end{enumerate}}
\newcommand{\bal}{\begin{align*}}
\newcommand{\eal}{\end{align*}}
\newcommand{\ba}{\begin{array}}
\newcommand{\ea}{\end{array}}
\newcommand{\eop}{$\square$\vspace{.3cm}}
\DeclareMathOperator{\card}{card}

\begin{document}
\title{The inverse problem for representation \\functions
of additive bases}
\author{Melvyn B. Nathanson\thanks{This work was supported
in part by grants from the NSA Mathematical Sciences Program
and the PSC-CUNY Research Award Program.}\\
Department of Mathematics\\
Lehman College (CUNY)\\
Bronx, New York 10468\\
Email: nathansn@alpha.lehman.cuny.edu}
\maketitle

\begin{abstract}
Let $A$ be a set of integers.
For every integer $n$, let $r_{A,2}(n)$ denote
the number of representations of $n$ in the form
$ n =  a_1 + a_2,$ where $a_1, a_2 \in A$ and $a_1 \leq a_2.$
The function $r_{A,2}: \Z \rightarrow \N_0\cup\{\infty\}$ 
is the {\em representation function of order $2$ for $A$.}
The set $A$ is called an {\em asymptotic basis of order $2$} if 
$r_{A,2}^{-1}(0)$ is finite, that is, if every integer 
with at most a finite number of exceptions can be represented 
as the sum of two not necessarily distinct elements of $A$.
It is proved that every function is a representation function, that is,
if $f: \Z \rightarrow \N_0\cup\{\infty\}$ is any function such that 
$f^{-1}(0)$ is finite, 
then there exists a set $A$ of integers such that $f(n) = r_{A,2}(n)$
for all $n \in \Z$.
Moreover, the set $A$ can be constructed so that 
$\card\{a\in A : |a| \leq x\} \gg x^{1/3}.$
\end{abstract}

\section{Representation functions}
Let \N, $\N_0,$ and \Z\ denote the positive integers, nonnegative integers, and 
integers, respectively.
Let $A$ and $B$ be sets of integers.  We define the {\em sumset}
\[
A+B = \{a+b : \text{$a\in A$ and $b \in B$} \},
\]
and, in particular, 
\[
2A = A+A = \{a_1 + a_2 : a_1, a_2 \in A \}
\]
and
\[
A+b = A+\{b\}= \{a+b : a \in A  \}.
\]
The {\em restricted sumsets} are
\[
A\hat{+} B = \{a+b : \text{$a\in A$, $b \in B$, and $a \neq b$} \}
\]
and
\[
2\wedge A = A \hat + A
= \{a_1 + a_2 : a_1, a_2 \in A\text{ and } a_1 \neq a_2 \}.
\]
Similarly, we define the {\em difference set}
\[
A-B = \{a-b : \text{$a\in A$ and $b \in B$} \}
\]
and 
\[
-A = \{0\}-A =  \{-a : -a \in A  \}.
\]
We introduce the {\em counting function}
\[
A(y,x) = \sum_{\substack{a\in A \\y \leq a \leq x}}1.
\]
Thus, $A(-x,x)$ counts the number of elements 
$a \in A$ such that $|a| \leq x.$

For functions $f$ and $g$, we write $f \gg g$ 
if there exist numbers $c_0$ and $x_0$
such that $|f(x)| \geq c_0|g(x)|$ for all $x \geq x_0,$ 
and $f \ll g$ if $|f(x)| \leq c_0|g(x)|$ for all $x \geq x_0.$ 

In this paper we study representation functions of sets of integers.
For any set $A \subseteq \Z,$
the {\em representation function} $r_{A,2}(n)$ counts the number of
ways to write $n$ in the form $n = a_1 + a_2$, 
where $a_1, a_2 \in A$
and $a_1 \leq a_2$.  
The set $A$ is called an {\em asymptotic basis of order 2}
if all but finitely many integers can be represented 
as the sum of two not necessarily distinct elements of $A$,
or, equivalently, if the function 
\[
r_{A,2}: \Z \rightarrow \N_0 \cup \{\infty\}
\]
satisfies 
\[
\card(r_{A,2}^{-1}(0)) < \infty.
\]

Similarly, the {\em restricted representation function} $\hat{r}_{A,2}(n)$ 
counts the number of ways to write $n$ in the form $n = a_1 + a_2$, 
where $a_1, a_2 \in A$ and $a_1 < a_2$.  
The set $A$ is called a {\em restricted asymptotic basis of order 2}
if all but finitely many integers can be represented 
as the sum of two distinct elements of $A$.

Let 
\beq           \label{invr:f1}
f: \Z \rightarrow \N_0 \cup \{\infty\}
\eeq
be any function 
such that 
\beq           \label{invr:f2}
\card(f^{-1}(0)) < \infty.
\eeq
The {\em inverse problem for representation functions of order $2$}
is to find sets $A$ such that $r_{A,2}(n) = f(n)$ for all $n \in \Z.$
Nathanson~\cite{nath03c} proved that  
every function $f$ satisfying~(\ref{invr:f1})
and~(\ref{invr:f2}) is the representation
function of an asymptotic basis of order 2,
and that such bases $A$ can be arbitrarily {\em thin}
in the sense that the counting functions
$A(-x,x)$ tend arbitrarily slowly to infinity.
It remained an open problem to construct {\em thick}
asymptotic bases of order 2 for the integers
with a prescribed representation function.

In the special case of the function $f(n) = 1$ for all integers $n$,
Nathanson~\cite{nath03a} constructed a unique representation basis, 
that is, a set $A$ of integers with $r_{A,2}(n) = 1$ for all $n\in \Z$, 
with the additional property that $A(-x,x) \gg \log x$. 
He posed the problem of constructing a unique representation basis $A$ 
such that $A(-x,x) \gg x^{\alpha}$ for some $\alpha > 0.$

In this paper we prove that for {\em every} function $f$ 
satisfying~(\ref{invr:f1}) and~(\ref{invr:f2}) 
there exist uncountably many asymptotic bases $A$ of order 2
such that $r_{A,2}(n) = f(n)$ for all $n \in \Z,$ and $A(-x,x) \gg x^{1/3}$.
It is not known if there exists a real number $\delta > 0$ such that
one can solve the inverse problem for arbitrary functions $f$
satisfying~(\ref{invr:f1}) and~(\ref{invr:f2}) 
with $A(-x,x) \gg x^{1/3 + \delta}$.

\section{The Erd\H os-Tur\' an conjecture}
The set $A$ of nonnegative integers is an asymptotic basis 
of order 2 for $\N_0$
if the sumset $2A$ contains all sufficently large integers.  
If $A$ is a set of nonnegative integers, then 
\[
0 \leq r_{A,2}(n) < \infty
\]
for every $n \in \N_0$.  It is not true, however, that if 
\[
f:\N_0 \rightarrow \N_0
\]
is a function with
\[
\card\left( f^{-1}(0) \right) < \infty,
\]
then there must exist a set $A$ of nonnegative integers such that 
$r_{A,2}(n) = f(n)$ for all $n \in \N_0$.
For example, Dirac~\cite{dira51} proved that the representation function
of an asymptotic basis of order 2 cannot be eventually constant,
and Erd\H os and Fuchs~\cite{erdo-fuch56} proved that the mean value
$\sum_{n \leq x}r_{A,2}(n)$ of an asymptotic basis of order 2 cannot
converge too rapidly to $c x$ for any $c > 0.$
A famous conjecture of Erd\H os and Tur\' an~\cite{erdo-tura41} 
states that the representation function of an asymptotic 
basis of order 2 must be unbounded.
This problem is only a special case of the general
{\em inverse problem for representation functions for bases 
for the nonnegative integers}:  Find necessary and sufficient conditions
for a function $f:\N_0 \rightarrow \N_0$ 
satisfying $\card\left( f^{-1}(0) \right) < \infty$
to be the representation function of an asymptotic basis of order $2$ for $\N_0$.

It is a remarkable recent discovery that 
the inverse problem for representation functions 
for the integers, and, more generally, 
for arbitrary countably infinite abelian groups
and countably infinite abelian semigroups with a group component, 
is significantly easier than 
the inverse problem for representation functions 
for the nonnegative integers and for other 
countably infinite abelian semigroups (Nathanson~\cite{nath04a}).

\section{Construction of thick bases for the integers}

Let $[x]$ denote the integer part of the real number $x$.

\bl           \label{invr:lemma:uk}
Let $f:\Z \rightarrow \N_0\cup\{\infty\}$ be a function such that $f^{-1}(0)$ is finite.  
Let $\Delta$ denote the cardinality of the set $f^{-1}(0).$
Then there exists a sequence $U = \{u_k\}_{k=1}^{\infty}$
of integers such that, for every $n \in \Z$ and $k \in \N$,
\[
f(n) = \card\{ k\geq 1 : u_k = n \}
\]
and
\[
|u_k| \leq \left[\frac{k+\Delta}{2}\right].
\]
\el

\pf
Every positive integer $m$ can be written uniquely in the form 
\[
m = s^2+s+1+r, 
\]
where $s$ is a nonnegative integer and $|r|\leq s.$
We construct the sequence 
\begin{align*}
V & = \{0,-1,0,1,-2,-1,0,1,2,-3,-2,-1,0,1,2,3,\ldots\} \\
& = \{v_m\}_{m=1}^{\infty},
\end{align*}
where 
\[
v_{s^2+s+1+r} = r \qquad\mbox{for $|r|\leq s.$}
\]
For every nonnegative integer $k$, the first occurrence of $-k$ 
in this sequence is $v_{k^2+1} =  -k,$
and the first occurrence of $k$ in this sequence is $v_{(k+1)^2} =  k.$

The sequence $U$ will be the unique subsequence of $V$ constructed as follows.
Let $n \in \Z.$  If $f(n) = \infty,$ then $U$ will contain 
the terms $v_{s^2+s+1+n}$ for every $s \geq |n|$.
If $f(n) = {\ell} <\infty,$ then $U$ will contain 
the ${\ell}$ terms $v_{s^2+s+1+n}$ for $s = |n|, |n|+1,\ldots,|n|+{\ell}-1$ in the subsequence $U$,
but not the terms $v_{s^2+s+1+n}$ for $s \geq |n| + {\ell}.$
Let $m_1 < m_2 < m_3 < \cdots$ be the strictly increasing sequence of positive integers
such that $\{v_{m_k}\}_{k=1}^{\infty}$ is the resulting subsequence of $V$.
Let $U = \{u_k\}_{k=1}^{\infty}$, where $u_k = v_{m_k}.$
Then 
\[
f(n) = \card\{ k \geq 1 : u_k = n \}.
\]

Let $\card\left(f^{-1}(0)\right) = \Delta.$
The sequence $U$ also has the following property:
If $|u_k| = n,$ then for every integer $m \not\in f^{-1}(0)$ 
with $|m| < n$ there is a positive integer $j < k$ with $u_j = m$.  
It follows that
\[
\{0,1,-1,2,-2,\ldots, n-1, -(n-1)\}\setminus f^{-1}(0) \subseteq \{u_1,u_2,\ldots,u_{k-1}\},
\]
and so
\[
k-1 \geq 2(n-1)+1 - \Delta.
\]
This implies that
\[
|u_k| = n \leq \frac{k+\Delta}{2}.
\]
Since $u_k$ is an integer, we have
\[
|u_k| \leq \left[\frac{k+\Delta}{2}\right].
\]
This completes the proof.
\eop

Lemma~\ref{invr:lemma:uk} is best possible in the sense that
for every nonnegative integer $\Delta$ there is a function
$f:\Z \rightarrow \N_0\cup\{\infty\}$ with $\card\left(f^{-1}(0)\right) = \Delta$ 
and a sequence $U = \{u_k\}_{k=1}^{\infty}$
of integers such that
\beq               \label{invr:U}
|u_k| = \left[\frac{k+\Delta}{2}\right]  \qquad \text{for all $k \geq 1.$}
\eeq
For example, if $\Delta = 2\delta + 1$ is odd, define the function $f$ by
\[
f(n) = \left\{\ba{ll}
0 & \text{if $|n| \leq \delta$}  \\
1 & \text{if $|n| \geq \delta + 1$} 
\ea
\right.
\]
and the sequence $U$ by
\begin{align*}
u_{2i-1} & = \delta + i, \\
u_{2i} & = -(\delta + i)
\end{align*}
for all $i \geq 1.$

If $\Delta = 2\delta$ is even, define $f$ by
\[
f(n) = \left\{\ba{ll}
0 & \text{if $-\delta \leq n \leq \delta-1$}  \\
1 & \text{if $n \geq \delta$ or $n \leq -\delta-1$} 
\ea
\right.
\]
and the sequence $U$ by $u_1 = \delta$ and
\begin{align*}
u_{2i} & = \delta + i, \\
u_{2i+1} & = -(\delta + i)
\end{align*}
for all $i \geq 1.$
In both cases the sequence $U$ satisfies~(\ref{invr:U}).

\bt                          \label{invr:theorem:thick}
Let $f:\Z \rightarrow \N_0 \cup \{\infty\}$
be any function such that 
\[
\Delta = \card(f^{-1}(0)) < \infty.
\]
Let
\[
c = 8 + \left[\frac{\Delta + 1}{2} \right].
\]
There exist uncountably many sets $A$ of integers such that
\[
r_{A,2}(n) = f(n) \qquad\text{for all $n \in \Z$}
\]
and
\[
A(-x,x) \geq \left(\frac{x}{c}\right)^{1/3}.
\]
\et

\pf
Let 
\[
\Delta = \card(f^{-1}(0)).
\]
By Lemma~\ref{invr:lemma:uk}, there exists a sequence 
$U = \{u_k\}_{k=1}^{\infty}$ of integers such that
\beq                         \label{invr:f}
f(n) = \card( \{i \in \N : u_i = n\})
\qquad\text{for all integers $n$}
\eeq
and
\beq                         \label{invr:u}
|u_k| \leq \frac{k+\Delta}{2} \qquad\text{for all $k \geq 1$.}
\eeq
We shall construct a strictly increasing sequence $\{i_k\}_{k=1}^{\infty}$
of positive integers and an increasing sequence $\{A_k\}_{k=1}^{\infty}$
of finite sets of integers such that, for all positive integers $k$, 
\benum
\item[(i)]
\[
|A_k| = 2k,
\]
\item[(ii)]
There exists a positive number $c$ such that 
\[
A_k \subseteq [-ck^3,ck^3]
\]
\item[(iii)]
\[
r_{A_k,2}(n) \leq f(n) \qquad\text{for all $n \in \Z$},
\]
\item[(iv)]
For $j = 1,\ldots, k,$
\[
r_{A_k,2}(u_j) \geq \card(\{i \leq i_k: u_i = u_j\}).
\]
\eenum

Let $\{A_k\}_{k=1}^{\infty}$ be a sequence of finite sets satisfying~(i)--(iv).
We form the infinite set
\[
A = \bigcup_{k=1}^{\infty} A_k.
\]
Let $x \geq 8c$, and let $k$ be the unique positive integer such that 
\[
ck^3 \leq x < c(k+1)^3.
\]
Conditions~(i) and~(ii) imply that 
\[
A(-x,x) \geq |A_k| = 2k > 2\left(\frac{x}{c}\right)^{1/3}-2 \geq \left(\frac{x}{c}\right)^{1/3}.
\]
Since
\[
f(n) = \lim_{k\rightarrow\infty} \card(\{i \leq i_k: u_i = n\}),
\]
conditions~(iii) and~(iv) imply that 
\[
r_{A,2}(n) = \lim_{k\rightarrow\infty}r_{A_k,2}(n) = f(n)
\]
for all $n \in \Z.$

We construct the sequence $\{A_k\}_{k=1}^{\infty}$ as follows.
Let $i_1 = 1.$  
The set $A_1$ will be of the form $A_1 = \{a_1+u_{i_1},-a_1\}$, 
where the integer $a_1$ is chosen so that $2A_1 \cap f^{-1}(0) = \emptyset$
and $a_1+u_{i_1} \neq -a_1$.
This is equivalent to requiring that
\beq                   \label{invr:a1exclude}
2a_1 \not\in (f^{-1}(0)-2u_{i_1}) \cup (-f^{-1}(0)) \cup \{-u_{i_1}\}.
\eeq
This condition excludes at most $1+2\Delta$ integers, 
and so we have at least two choices for the number $a_1$ 
such that $|a_1| \leq 1+\Delta$ and $a_1$ satisfies~(\ref{invr:a1exclude}).
Since $|u_{i_1}| = |u_1| \leq (1+\Delta)/2$ and 
\[
|a_1+u_{i_1}| \leq |a_1| + |u_{i_1}| \leq \frac{3(1 + \Delta)}{2},
\]
it follows that $A_1 \subseteq [-c,c]$ for any $c \geq 3(1+\Delta)/2,$
and the set $A_1$ satisfies conditions~(i)--(iv).

Let $k \geq 2$ and suppose that we have constructed sets
$A_1,\ldots,A_{k-1}$ and integers $i_1 < \cdots < i_{k-1}$
that satisfy conditions~(i)--(iv).
Let $i_k > i_{k-1}$ be the least integer such that 
\[
r_{A_{k-1},2}(u_{i_k}) < f(u_{i_k}).
\]
Since
\begin{align*}
i_k -1 & \leq \sum_{n\in \{ u_1,u_2,\ldots,u_{i_k-1} \}} r_{A_{k-1},2}(n)\\
& \leq \sum_{n\in \Z} r_{A_{k-1},2}(n) \\
& = \binom{2k-1}{2} \\
& < 2k^2,
\end{align*}
it follows that
\[
i_k \leq 2k^2.
\]
Also,~(\ref{invr:u}) implies that
\beq                   \label{invr:uik}
|u_{i_k}| \leq \frac{i_k+\Delta}{2} \leq k^2 + \frac{\Delta}{2}.
\eeq
We want to choose an integer $a_k$ such that the set
\[
A_k = A_{k-1} \cup \{ a_k+u_{i_k},-a_k \}
\]
satisfies~(i)--(iv).  
We have $|A_k| = 2k$ if
\[
a_k+u_{i_k} \neq -a_k
\]
and
\[
A_{k-1} \cap \{ a_k+u_{i_k},-a_k \} = \emptyset,
\]
or, equivalently, if
\beq                   \label{invr:akexclude}
a_k \not\in (-A_{k-1}) \cup (A_{k-1} - u_{i_k}) \cup \{-u_{i_k}/2\}.
\eeq
Thus, in order for $A_{k-1} \cup \{ a_k+u_{i_k},-a_k \}$ to satisfy
condition~(i), we exclude at most $2|A_{k-1}|+1 = 4k-3$ integers
as possible choices for $a_k.$

The set $A_k$ will satisfy conditions~(iii) and~(iv) if
\[
2A_k \cap f^{-1}(0) = \emptyset
\]
and
\[
r_{A_k,2}(n) = 
\left\{
\begin{array}{ll}
r_{A_{k-1},2}(n) & \text{for all $n \in 2A_{k-1}\setminus\{u_{i_k}\}$} \\
r_{A_{k-1},2}(n)+1 & \text{for $n = u_{i_k}$} \\
1 & \text{for all $n \in 2A_k \setminus \left(2A_{k-1}\cup \{u_{i_k}\} \right).$}
\end{array}
\right.
\]
Since the sumset $2A_k$ decomposes into
\begin{align*}
2A_k& = 2\left( A_{k-1} \cup \{ a_k+u_{i_k},-a_k \} \right) \\
& = 2A_{k-1} \cup \left( A_{k-1} + \{ a_k+u_{i_k},-a_k \} \right)
\cup \{u_{i_k}, 2a_k+2u_{i_k}, -2a_k\},
\end{align*}
it suffices that
\begin{align}
\left( A_{k-1} + \{ a_k+u_{i_k},-a_k \} \right) 
\cap 2A_{k-1} & = \emptyset,    \label{invr:exc1}        \\
\left( A_{k-1} + \{ a_k+u_{i_k},-a_k \} \right) 
\cap f^{-1}(0) & = \emptyset,     \label{invr:ecl2}       \\
\left( A_{k-1} + a_k+u_{i_k} \right) 
\cap \left( A_{k-1} -a_k \right) 
& = \emptyset,  \label{invr:ecl3}       \\
\{2a_k+2u_{i_k},-2a_k \} \cap 2A_{k-1} & = \emptyset   \label{invr:ecl4}\\
\{2a_k+2u_{i_k},-2a_k \} \cap  f^{-1}(0) & = \emptyset   \label{invr:ecl5}\\
\{2a_k+2u_{i_k},-2a_k \} \cap \left(A_{k-1} + \{ a_k+u_{i_k},-a_k \} \right)
 & = \emptyset.   \label{invr:ecl6}
\end{align}
Equation~(\ref{invr:exc1}) implies that the integer $a_k$ must be chosen so
that it cannot be represented either in the form 
\[
a_k = x_1 + x_2 - x_3 - u_{i_k}
\]
or
\[
a_k = x_1 - x_2 - x_3,
\]
where $x_1,x_2,x_3 \in A_{k-1}$.
Since $\card(A_{k-1}) = 2(k-1),$ it follows that the number of integers 
that cannot be chosen as the integer $a_k$ because of 
equation~(\ref{invr:exc1}) is at most $2(2(k-1))^3 = 16(k-1)^3$.

Similarly, the numbers of integers excluded as possible choices for $a_k$ because of 
equations~(\ref{invr:ecl2}),~(\ref{invr:ecl3}),~(\ref{invr:ecl4}),~(\ref{invr:ecl5}), 
and~(\ref{invr:ecl6}) are at most $4\Delta(k-1), 4(k-1)^2, 8(k-1)^2, 2\Delta,$ 
and $8(k-1),$ respectively, and so the number of integers 
that cannot be chosen as $a_k$ is 
\begin{align*}
16(k-& 1)^3 + 12(k-1)^2 + (4\Delta + 8)(k-1) + 2\Delta \\
& = 16k^3 - 36k^2 + (32+4\Delta)k - 2\Delta - 12\\
& \leq (16+\Delta)k^3 - 4k^2 -32k(k-1)- 2\Delta -12.
\end{align*}
Let
\[
c = 8 + \left[ \frac{\Delta + 1}{2} \right].
\]
The number of integers $a$ with
\beq           \label{invr:a}
|a| \leq ck^3 - k^2 -  \left[ \frac{\Delta + 1}{2} \right] =
\left( 8 + \left[ \frac{\Delta + 1}{2} \right] \right) k^3
 - k^2 -  \left[ \frac{\Delta + 1}{2} \right] 
\eeq
is 
\begin{align*}
\left( 16 +  2 \left[\frac{\Delta + 1}{2} \right] \right) k^3 
& - 2k^2 -  2\left[ \frac{\Delta + 1}{2} \right]  + 1 \\
& \geq \left( 16 + \Delta \right) k^3 - 2k^2 -  \Delta.  
\end{align*}
If the integer $a$ satisfies~(\ref{invr:a}), then~(\ref{invr:uik}) implies that
\[
|a+u_{i_k}| \leq |a| +| u_{i_k}| \leq ck^3.
\]
It follows that there are at least two acceptable choices of the integer $a_k$
such that the set $A_k = A_{k-1} \cup \{ a_k+u_{i_k}, -a_k \}$ 
satisfies conditions~(i)--(iv).
Since this is true at each step of the induction, there are uncountably many sequences
$\{A_k\}_{k=1}^{\infty}$ that satisfy conditions~(i)--(iv).
This completes the proof.
\eop

We can modify the proof of Theorem~\ref{invr:theorem:thick}
to obtain the analogous result for the restricted representation
function $\hat{r}_{A,2}(n).$

\bt                          \label{invr:theorem:thick2}
Let $f:\Z \rightarrow \N_0 \cup \{\infty\}$
be any function such that 
\[
\card(f^{-1}(0)) < \infty.
\]
Then there exist uncountably many sets $A$ of integers such that
\[
\hat{r}_{A,2}(n) = f(n) \qquad\text{for all $n \in \Z$}
\]
and
\[
A(-x,x) \gg x^{1/3}.
\]
\et

\section{Representation functions for bases of order $h$}

We can also prove similar results for the representation functions 
of asymptotic bases and restricted asymptotic bases of order $h$ 
for all $h \geq 2$.

For any set $A \subseteq \Z,$
the {\em representation function} $r_{A,h}(n)$ counts the number of
ways to write $n$ in the form $n = a_1 + a_2 + \cdots + a_h$, 
where $a_1, a_2, \ldots, a_h \in A$
and $a_1 \leq a_2 \leq \cdots \leq a_h$.  
The set $A$ is called an {\em asymptotic basis of order h}
if all but finitely many integers can be represented 
as the sum of $h$ not necessarily distinct elements of $A$,
or, equivalently, if the function 
\[
r_{A,h}: \Z \rightarrow \N_0 \cup \{\infty\}
\]
satisfies 
\[
\card(r_{A,h}^{-1}(0)) < \infty.
\]
Similarly, the {\em restricted representation function} $\hat{r}_{A,h}(n)$ 
counts the number of ways to write $n$ as a sum of $h$ pairwise distinct elements of $A$.
The set $A$ is called a {\em restricted asymptotic basis of order h}
if all but finitely many integers can be represented 
as the sum of $h$ pairwise distinct elements of $A$.

\bt                          \label{invr:theorem:thick-h}
Let $f:\Z \rightarrow \N_0 \cup \{\infty\}$
be any function such that 
\[
\card(f^{-1}(0)) < \infty.
\]
There exist uncountably many sets $A$ of integers such that
\[
r_{A,h}(n) = f(n) \qquad\text{for all $n \in \Z$}
\]
and
\[
A(-x,x) \gg x^{1/(2h-1)},
\]
and there exist uncountably many sets $A$ of integers such that
\[
\hat{r}_{A,h}(n) = f(n) \qquad\text{for all $n \in \Z$}
\]
and
\[
A(-x,x) \gg x^{1/(2h-1)}.
\]
\et

\providecommand{\bysame}{\leavevmode\hbox to3em{\hrulefill}\thinspace}
\providecommand{\MR}{\relax\ifhmode\unskip\space\fi MR }
\providecommand{\MRhref}[2]{%
  \href{http://www.ams.org/mathscinet-getitem?mr=#1}{#2}
}
\providecommand{\href}[2]{#2}

\end{document}